\pdfoutput=1

\documentclass[onecolumn]{myautart}    

\usepackage{graphicx}          
\usepackage{color}

\usepackage{amsmath}

\usepackage{amssymb}  

\usepackage{algpseudocode} 

\newcommand{\bbaru}{\boldsymbol{\bar{u}}}
\newcommand{\bbarx}{\boldsymbol{\bar{x}}}

\newcommand{\bu}{\boldsymbol{u}}

\newcommand{\bx}{\boldsymbol{x}}

\usepackage{pgfplots}
\newlength\figureheight
\newlength\figurewidth

\newtheorem{problem}{Problem}
\DeclareMathOperator{\atan}{atan}

\DeclareMathOperator{\Id}{Id}

\begin{document}

\begin{frontmatter}

\title{Finite-dimensional control of linear discrete-time \\fractional-order systems}

\author[IST]{Andrea Alessandretti}\ead{andrea.alessandretti@fe.up.pt}\;\;\;
\author[REN]{S\'{e}rgio Pequito}\ead{goncas@rpi.edu}\;\;\;
\author[PEN]{George J. Pappas}\ead{pappasg@seas.upenn.edu}\;\;\;
\author[IST]{A. Pedro Aguiar}\ead{pedro.aguiar@fe.up.pt} 

\address[IST]{Faculty of Engineering, University of Porto, Porto, Portugal.}      
\address[REN]{Department of Industrial and Systems Engineering, Rensselaer Polytechnic Institute,Troy, USA}       
\address[PEN]{Department of Electrical \& Systems Engineering School, Engineering and Applied Science, University of Pennsylvania, USA}       



\begin{abstract}                          
This paper addresses the design of \mbox{finite-dimensional} feedback control laws for linear discrete-time fractional-order systems with additive state disturbance. A set of sufficient conditions are provided to guarantee convergence of the state trajectories to an ultimate bound around the origin with size increasing with the magnitude of the disturbances. Performing a suitable change of coordinates, the latter result can be used to design a controller that is able to track reference trajectories that are solutions of the unperturbed \mbox{fractional-order} system. To overcome the challenges associated with the generation of such solutions, we address the practical case where the references to be tracked are generated as a solution of a specific \mbox{finite-dimensional} approximation of the original fractional-order system. In this case, the tracking error trajectory is driven to an asymptotic bound that is increasing with the magnitude of the disturbances and decreases with the increment in the accuracy of the approximation. The proposed controllers are finite-dimensional, in the sense that the computation of the control input only requires a finite number of previous state and input vectors of the system. Numerical simulations illustrate the proposed design methods in different scenarios.
\end{abstract}

\end{frontmatter}

\section{Introduction}

A fractional-order derivative is a generalization of an integer order derivative, where the order can take non-integer values. Its discretization results in a fractional-order difference \cite{samko1993fractional} that rules the dynamics of the discrete time counterpart. In this paper, we address the design of control laws for dynamical systems described by fractional-order difference equations, i.e., discrete-time fractional-order systems (FOSs).

In contrast to integer order differences, which use local information around a given point, fractional-order differences are obtained using all the previous values of a given function up to the point where the difference is computed. Because of the global nature of the fractional-order operators, models described by fractional-order differences seem to be particularly suitable to capture the behavior of systems with infinite memory and hence have found applications in many fields. Many biological systems manifest fractional-order behaviors \cite{Pereira2010}. \mbox{Fractional-order} derivatives are used to model the cardiac \mbox{tissue-electrode} interface and electrically stimulated or mechanically stressed tissues \cite{Magin2008,Magin2010}, drugs relies and adsorption \cite{Magin2011}, and the behavior of neural systems \cite{Anastasio1994}. Applications on polymer science are highlighted in \cite{Hilfer2000}. In \cite{Sun1984}, transfer functions with a pole of fractional power are used to model analog electronic circuits. The work \cite{Cao2010} explores the benefits of fractional models for modeling of fuel cells.

The many applications and the control challenges associated with fractional-order systems motivated a wide literature in the control community. For the \mbox{continuous-time} case, \cite{Matignon96stabilityresults} shows the equivalence of stability with the presence of the eigenvalues of a given polynomial in a specific angular sector, generalizing the integer order case where such a sector becomes the left half plane. Building on this, the works in \cite{Ahn2008,Farges2010,NDoye2013} focus on developing linear matrix inequality (LMI) conditions for stability and extensions to the case of uncertain polytonic FOSs. The work \cite{739144} provides analytical expressions for the unit-step and unit-impulse of FOSs, that is then used to compare a proposed \mbox{fractional-order} PID controller with the traditional integer-order counterpart. The work \cite{Li2009a} presents the concept of Mittag-Leffler stability for nonlinear FOSs and proposes a Lyapunov-like condition to certify such stability, using both integer and fractional-order derivative of the Lyapunov function. 

For the discrete-time FOSs, the control design problem is addressed either using the infinite-dimensional formulation of the FOS \cite{Guermah2010}, or adopting a \mbox{finite-dimensional} approximation of the FOS obtained by truncating the infinite-dimensional \mbox{fractional-order} difference operator and using the associated stabilizing controller to control the original FOS \cite{kaczorek2008practical,Sopasakis2017}. The main drawback of the first approach is that the stability conditions that arise from the analysis of the whole infinite-dimensional form of FOS are generally difficult to verify and to apply for (finite-dimensional) control design purposes. The second approach was initially introduced in \cite{kaczorek2008practical},  without providing stability guarantees, and then adopted in the work \cite{Sopasakis2017} for the design of a stabilizing Model Predictive Control scheme. In the latter approach, the authors derive sufficient conditions for practical asymptotic stability and, upon satisfaction of an additional condition, asymptotic stability of the origin. In addition, the latter work is able to explicitly handle constraints in both the system state and input vectors.

In this work, we retain the simplicity in the design of the widely used approach proposed in \cite{kaczorek2008practical}, while at the same time providing closed-loop guarantees. The analysis adopted in this work builds upon the reasoning proposed in \cite{Sopasakis2017}. Notwithstanding, as the main differences, this work does not consider system constraints and shows that convergence can be guaranteed globally, i.e., for any initial conditions of the system. The proposed controller is easy to implement, without the need of optimization-based schemes. In addition, we address (\emph{i}) possible additive disturbances acting on the system and (\emph{ii}) the trajectory-tracking problem. For the latter case, we also consider the practical scenario where the reference trajectories are solutions of an exogenous finite-dimensional approximation of the original infinite-dimensional FOS. Specifically, we provide sufficient conditions for the tracking error to converge to an asymptotic bound that is increasing with the magnitude of the disturbances and decreases with the increment in the accuracy of the approximation.

The remainder of this paper is structured as follows. Section~\ref{s:prob} introduces the linear discrete-time fractional-order system and formalizes the control problem addressed in this paper. Section~\ref{s:main} contains the main results that are illustrated in Section~\ref{s:numeric} with numerical simulations. All the proofs are relegated to the Appendix.

{\bf Notation:} The terms $\mathbb{R}$, $\mathbb{R}^+$, $\mathbb{Z}$, $\mathbb{N}$, $\mathbb{N}^{+}$ denote the set of reals, non-negative reals, integers, non-negative integers, and positive integers numbers, respectively. The terms $\mathbb{R}^n$ and $\mathbb{R}^{n\times m}$ represent the set of column vectors of size $n$ and $n$-by-$m$ matrices with real entries. A continuous function $\alpha:\mathbb{R}^+\to \mathbb{R}^+$ is said to belong to class $\mathcal{K}$ if it is strictly increasing and $\alpha(0)=0$. In addition, a function is said to belong to class $\mathcal{K}_\infty$ if, in addition, $\alpha(r)\to+\infty$ as $r \to +\infty$. Moreover, a continuous function $\beta:\mathbb{R}^+\times \mathbb{R}^+\to \mathbb{R}^+$ is said to belong to class $\mathcal{KL}$ if, for a constant~$s$, $\beta(r,s)$ belongs to class $\mathcal{K}$ with respect to $r$, and for a constant $r$ it is decreasing on $s$ and satisfies $\beta(r,s)\to 0$ as $s\to\infty$ (see \cite{khalil2002nonlinear}). The term $\Id$ denotes the identity function, i.e., $\Id(x)=x$. The term $I_n$ is used to denoted a $n\times n$ identity matrix. For a sequence $\{x(k)\}_{k\in\mathbb{Z}}$, we use the notation $\|x\|_\infty=\sup_{k\in\mathbb{Z}}\|x(k)\|$ and $\|x\|_{[k_1,k_2]}=\max_{k_1\leq k \leq k_2}\|x(k)\|$, if $k_1,k_2\in\mathbb{Z}$ with $k_1\leq k_2$, and $\|x\|_{[k_1,k_2]}=0$ if $k_1> k_2$. Given a square matrix $M\in\mathbb{R}^{n\times n}$, the notation $M\succ0$ indicates that the matrix is positive definite, i.e., $v^\top M v>0$ for any $v\in\mathbb{R}^n$. Moreover, whenever clear from the context, we use $x^+$ and $x$ to refer to $x(k+1)$ and $x(k)$, respectively.

\section{Control problem formulation}\label{s:prob}

This section introduces the discrete-time FOS and defines the control objectives addressed in this paper.

\subsection{Discrete-time fractional-order systems}
Consider a left-bounded sequence $\{z(k)\}_{k\in\mathbb{Z}}$ over~$k$, i.e., $\limsup\limits_{k\to -\infty}\|z(k)\|<\infty$. Then, for any ${{a}}\in\mathbb{R}^+$ the Gr\"{u}nwald-Letnikov \mbox{fractional-order} difference \cite{samko1993fractional} is defined as \vspace{-0.25cm}
\begin{align}
\Delta^{{{a}}}z(k):=\sum_{j=0}^{\infty}c_j^{{{a}}} z(k-j), && c_j^{{{a}}} = (-1)^j\begin{pmatrix}{{a}}\\j\end{pmatrix},\label{eq:defDelta}
\end{align}
where \vspace{-0.25cm}
\begin{align}
\begin{pmatrix}{{a}}\\j\end{pmatrix}=
\begin{cases}1,& j=0\\
\prod_{i=0}^{j-1} {{{a}}-i \over i+1}& j>0\end{cases}
\end{align}
for all $j\in\mathbb{N}$. Notice that the summation in \eqref{eq:defDelta} is well defined from the uniform boundedness of the sequence $\{z(k)\}_{k\in\mathbb{Z}}$ and the fact that \vspace{-0.25cm}
\begin{align}
|c_j^{{{a}}}|\leq{{{a}}^j\over j!},\label{eq:boundc}
\end{align} which implies that the sequence $\{c_j^{{{a}}}\}_{j\in\mathbb{N}}$ is absolutely summable for any ${{a}} \in \mathbb{R}^+$. A \emph{discrete-time fractional-order system with additive disturbance} can be modeled as  \vspace{-0.25cm}
\begin{align}
\sum_{i=1}^{l}A_i\Delta^{{{a}}_i}x(k+1) = \sum_{i=1}^{r}B_i\Delta^{{{b}}_i}u(k)+  \sum_{i=1}^{s}G_i\Delta^{{{g}}_i}w(k)\label{eq:sys}
\end{align}
for scalars ${{a}}_i\in \mathbb{R}^+$, with $i=1,\dots,l$, ${{b}}_i\in \mathbb{R}^+$, with $i=1,\dots,r$, and ${{g}}_i\in \mathbb{R}^+$, with $i=1,\dots,s$, and where $x(k)\in\mathbb{R}^n$, $u(k)\in\mathbb{R}^m$, and $w(k)\in\mathbb{R}^{p}$ denote the state, input, and disturbance vectors at time~$k\in \mathbb{N}$, respectively. The latter vector is considered to be bounded as \begin{align}
\|w(k)\|\leq b_w, &&k\in\mathbb{N},\label{eq:bound_w}
\end{align} for some positive scalar $b_w$. We denote by $x_0=x(0)$ the initial condition of the FOS at time~$k=0$. In the computation of the fractional-order difference, the state, input, and disturbances are considered to be zero before the initial time, i.e., $x(k)=0$, $u(k)=0$, and $w(k)=0$ for all $k<0$. 

It is worth noticing that the infinite dimensional case where the state is defined for all the times, can be addressed in the presented theorems considering the limit of the initial time that goes to minus infinity (or equivalently, considering the case $k\to +\infty$). In addition, given the robustness properties of the proposed control law, it is possible to addressed the case where the controller is activated at a given time by considering the previous absence of control input (or zero input signal) as a (bounded) disturbance. 

\subsection{Control objectives}
The first control objective concerns the regulation of the state of the FOS toward the origin.
\begin{problem}[Convergence]\label{pb:1} Design a \linebreak \mbox{finite-dimensional} controller that renders the associated closed-loop system with \eqref{eq:sys} globally bounded and globally ultimately bounded with respect to the disturbance \eqref{eq:bound_w}, i.e., for any initial condition $x_0\in\mathbb{R}^n$ the associated closed-loop state trajectory is bounded and satisfies
\begin{align}
\|x(k)\|\leq\beta(\|x_0\|,k) + \gamma(b_w), && k\in\mathbb{N} \label{eq:issdef}
\end{align}
for a function $\beta:\mathbb{R}^+\times \mathbb{R}^+\to \mathbb{R}^+$, with $\beta(r,s)\to 0$ as $s\to +\infty$, and a class-$\mathcal{K}$ function $\gamma:\mathbb{R}^+\to \mathbb{R}^+$. 
\end{problem}

\begin{defn}[Finite-dimensional controller]
We denote by \emph{finite-dimensional controller} a control law that only requires a finite number of previous state and input vectors to decide the current input to apply to the system, i.e., there exist two values $v_1,v_2\in \mathbb{N}$ such that the control law can be expressed in the form $u(k)=\kappa(x(k),\dots,x(k-v_1),u(k-1),\dots u(k-v_2))$.
\end{defn}

As it will be made explicit in the remainder of the paper, performing a suitable change of coordinate, a controller that solves Problem \ref{pb:1} can be used to track any reference trajectory that is a solution of \eqref{eq:sys} (similar to integer-order linear systems). It is interesting to note here that since \eqref{eq:sys} is infinite-dimensional, it often presents a major challenge in cases where one wishes to generate such a reference trajectory. Specifically, for any desired input reference trajectory~$\{u_r(k)\}_{k\in\mathbb{N}}$, the computation of the associated (nominal) fractional-order state trajectory $\{x_r(k)\}_{k\in\mathbb{N}}$ (i.e., the solution of \eqref{eq:sys} considered with $w(k)=0$, $u(k)=u_r(k)$ and $x(k)=x_r(k)$ for all $k\in\mathbb{N}$) becomes intractable as $k\to\infty$ due to the sum of an always growing number of terms. To avoid this situation, we consider the case where the reference trajectory is generated using a finite-dimensional exogenous system.

\begin{problem}[Tracking exogenous solutions]\label{pb:2} Let $x_e(k)\in\mathbb{R}^{n_e}$, $u_r(k)\in\mathbb{R}^m$, and $x_r(k)\in\mathbb{R}^{n}$, with $k\in\mathbb{N}$, be the state, input, and output vectors, respectively, of the \mbox{finite-dimensional} exogenous system
\begin{subequations}\label{eq:exosys}
\begin{align}
x_e(k+1) = Ax_e(k)+Bu_r(k),&& x_e(0)=x_{e,0},\\
x_r(k) = Cx_e(k),
\end{align}
\end{subequations}
where $x_{e,0}\in\mathbb{R}^{n_e}$ denotes the initial condition at time~$k=0$. Design the matrices $A$, $B$, and $C$ of \eqref{eq:exosys} and a \mbox{finite-dimensional} controller for the FOS \eqref{eq:sys} such that, for the associated closed-loop, the tracking error 
\begin{align}
e(k):=x(k)-x_r(k),
\end{align}
defined as the difference between the state of the FOS and the reference signal, satisfies 
\begin{align}
\|e(k)\|\leq\beta(\|e(0)\|,k) + \gamma(b_w) + d \label{eq:iips}
\end{align}
for a function $\beta:\mathbb{R}^+\times \mathbb{R}^+\to \mathbb{R}^+$ , with $\beta(r,s)\to 0$ as $s\to +\infty$, a class-$\mathcal{K}$ function $\gamma:\mathbb{R}^+\to \mathbb{R}^+$, and a constant $d>0$. 
\end{problem}

It is worth noticing that although the desired controller is finite-dimensional, the closed-loop guarantees of Problem~\ref{pb:2} hold when the plant is described by the closed-loop with the infinite-dimensional system \eqref{eq:sys}. 

 \section{Main contribution}\label{s:main}
Before stating the main results, different formulations of the FOSs \eqref{eq:sys} are introduced, as well as mild technical assumptions.

\begin{assum}\label{ass:model} The matrix $\sum_{i=1}^lA_i$ is invertible.
\end{assum}

Using Assumption \ref{ass:model}, the FOS model \eqref{eq:sys} can be reformulated as  \vspace{-0.25cm}\begin{align}
x(k+1) &= \sum_{j=1}^{\infty}\check{A}_jx(k-j+1)+\sum_{j=0}^{\infty}\check{B}_ju(k-j) +\sum_{j=0}^{\infty}\check{G}_jw(k-j)\label{eq:sys2}
\end{align}
where $\check{A}_j = -\hat{A}_{0}^{-1}\hat{A}_{j}$, $\check{B}_j = \hat{A}_{0}^{-1}\hat{B}_{j}$, and $\check{G}_j = \hat{A}_{0}^{-1}\hat{G}_{j}$ with $\hat{A}_{j}=\sum_{i=1}^lA_ic_j^{{{a}}_i}$, $\hat{B}_{j}=\sum_{i=1}^rB_ic_j^{{{b}}_i}$, and $\hat{G}_{j}~=~\sum_{i=1}^sG_ic_j^{{{g}}_i}$. 

Formulation \eqref{eq:sys2} highlights one of the main peculiarities of FOSs, i.e., that the state $x(k+1)$ does not depend only on the state, input, and disturbance at time $k$, but it is a function of the whole past trajectory.


For any given integer $v\in\mathbb{N}^{+}$, the model \eqref{eq:sys} can be equivalently rewritten as 
\begin{subequations} \label{eq:augsys}
\begin{align}
\tilde{x}(k+1) &= \tilde{A}_v \tilde{x}(k)  + \tilde{B}_v u(k) +  \tilde{G}_v r(k), \quad\tilde{x}(0)=\tilde{x}_0   \label{eq:sysLin}\\
r(k)&=\sum_{j=v+1}^{\infty}\check{A}_jx(k-j+1)+\sum_{j=v+1}^{\infty}\check{B}_ju(k-j) +\sum_{j=0}^{\infty}\check{G}_jw(k-j)\label{eq:defd}
\end{align}
\end{subequations}
using the augmented state vector  $$\tilde{x}(k) = [x^\top (k),\dots,x^\top (k-v+1),u^\top (k-1),\dots,u^\top (k-v)]^\top \in\mathbb{R}^{v(n+m)}$$ and appropriate matrices $\tilde{A}_v$, $\tilde{B}_v$, and $\tilde{G}_v$, where $\tilde{x}_0= [x_0^\top ,0,\dots,0]^\top $ denotes the initial condition. Specifically, the terms $\{\check{A}_j\}_{j=1,\dots,v}$ and $\{\check{B}_j\}_{j=1,\dots,v}$ are used to form the matrices $\tilde{A}_v$ and $\tilde{B}_v$, whereas the remaining terms $\{\check{G}_j\}_{j=1,\dots,\infty}$ and the states and input components not included in $\tilde{x}(k)$ are collected in the term $\tilde{G}_v r(k)$.

\begin{defn}[v-approximation] \label{def:vApp} The v-approximation of the discrete-time FOS \eqref{eq:sys} is defined to be the system \eqref{eq:sysLin} but with $r(k)=0$ for all $k\in\mathbb{N}$.
\end{defn}

\begin{assum}[v-stabilizability]\label{ass:controllability} The pair $(\tilde{A}_v,\tilde{B}_v)$ of the v-approximation of the FOS system \eqref{eq:sys} is stabilizable. \end{assum}

Assumption \ref{ass:controllability} implies the existence of a matrix $K_v$ such that for any positive definite matrix $Q$ the equality
\begin{align}
A_K^\top PA_K - P + Q = 0,&&A_K=\tilde{A}_v+\tilde{B}_vK_v\label{eq:lyap}
\end{align}
holds for a positive definite matrix $P$. At this point, we are ready to state one of the main results of this paper.

\begin{thm} \label{th:main}  Consider the FOS system \eqref{eq:sys}, the reformulation \eqref{eq:augsys}, and let Assumptions \ref{ass:model} and \ref{ass:controllability} hold. Moreover, consider $v\in\mathbb{N}$ such that  \vspace{-0.1cm}
\begin{align}
&c_\Psi \Psi(v)<1 \label{eq:conditionEps}
\end{align}\vspace{-0.65cm}\\
with
\begin{align}
\Psi(v) \!=\!\!\!\sum_{i=1}^l\|\hat{A}_{0}^{-1}A_i\| \phi_{{{a}}_i}(v) \!+\!\! \sum_{i=1}^r \|\hat{A}_{0}^{-1}B_iK_v\|\phi_{{{b}}_i}(v)\label{eq:defPsi}
\end{align}
and $\phi_{{a}}(v)~=~e^{{a}} -\sum_{j=0}^{v} {{{a}}^j\over j!}$, $c_\Psi~=~\sqrt{c_2\over \hat{c}_4 c_\rho\lambda_{\mathrm{min}}(P)}$, 
$c_2~=~\lambda_{\mathrm{max}}(\tilde{G}_v^\top P \tilde{G}_v) +{\|\tilde{G}_v^\top PA_K\|^2\over \theta \lambda_{\mathrm{min}}(Q)} $, $c_4={(1-\theta)\lambda_{\mathrm{min}}(Q) \over \lambda_{\mathrm{max}}(P)}$, $\hat{c}_4=\min(c_4,\hat{\theta})$, for any $\theta,\hat{\theta},c_\rho\in(0,1)$.
Then, for any positive definite matrix $Q\succ 0$, there exists a pair of matrices $K_v$ and $P\succ 0$ such that \eqref{eq:lyap} holds and the controller \begin{align}
u(k)=K_v\tilde{x}(k), \label{eq:input}
\end{align}
with $\tilde{x}(k)=[x^\top (k),\dots,x^\top (k-v+1),u^\top (k-1),\dots,u^\top (k-v)]^\top$, solves Problem \ref{pb:1} with
$
\gamma(r) = c_\gamma r
$, where $\kappa \in(c_\Psi \Psi(v),1)$ and $c_\gamma~=~c_{\Psi}{\kappa \over 1-\kappa} \sum_{i=1}^s \|\hat{A}_{0}^{-1}G_i\|e^{{{g}}_i}$.
\end{thm}

Theorem \ref{th:main} provides a set of sufficient conditions that guarantee the state vector of the closed-loop system \eqref{eq:sys} with \eqref{eq:input} to converge to an ultimate bound with size proportional to the maximum magnitude of the noise. The controller \eqref{eq:input} is designed to be any linear controller that stabilizes the $v$-approximation of \eqref{eq:sys}, where $v$ satisfies condition \eqref{eq:conditionEps}. It is worth noticing that since the state of the $v$-approximation (with $v>1$) stores the previous state and input signals of the original fractional-order system, the control \eqref{eq:input} is not a standard state-feedback controller for the FOS, but rather uses feedback of a window of past states and inputs, where the window size increases with $v$.

In the control design for integer-order linear system, a controller designed to stabilize the origin of the state space can be used to stabilize the system around any of its solutions, i.e., performing a suitable change of coordinates, it is possible to stabilize the origin of the error space defined as the difference between the current state of the system and a solution of the system that we wish to track.

Similarly, performing a suitable change of coordinates the controller of Theorem \ref{th:main} can be used to track any feasible solution of the FOS \eqref{eq:sys}. This fact is made explicit in the following corollary.

\begin{cor}[Tracking FOS solution]\label{cor:1} Let the assumptions of Theorem \ref{th:main} hold and consider a generic state and input reference trajectory pair $\{x_r(k)\}_{k\in\mathbb{N}}$ and $\{u_r(k)\}_{k\in\mathbb{N}}$, respectively, that are the solution of \eqref{eq:sys} with $w(k)=0,$ for all $k\geq 0$, and the associated augmented state vector $\tilde{x}_r(k)= [x_r^\top (k),\dots,x_r^\top (k-v+1),u_r^\top (k-1),\dots,u_r^\top (k-v)]^\top $. 
Then, the closed-loop \eqref{eq:sys} with  
\begin{align}
u(k)=u_r(k)+K_v\tilde{e}(k),&& \tilde{e}(k) = \tilde{x}(k)-\tilde{x}_r(k)
\end{align}
leads to an error trajectory $e(k):=x(k)-x_r(k)$ that is globally bounded and globally ultimately bounded with respect to the disturbance \eqref{eq:bound_w}. Specifically, for any initial condition $e_0\in\mathbb{R}^n$ the error trajectory is bounded over time and satisfies
\begin{align}
\|e(k)\|\leq\beta(\|e_0\|,k) + \gamma(b_w), && k\in\mathbb{N} \label{eq:issdef}
\end{align}
for a function $\beta:\mathbb{R}^+\times \mathbb{R}^+\to \mathbb{R}^+$, with $\beta(r,s)\to 0$ as $s\to +\infty$, and a class-$\mathcal{K}$ function $\gamma:\mathbb{R}^+\to \mathbb{R}^+$. 
\end{cor}

Corollary \ref{cor:1} highlights the robustness of the proposed control strategy, where the origin of the error space is ultimately bounded with the size of the ultimate bound being increasing with the bound of the disturbance.

In what follows, we consider the case where the trajectory that we want to track is not a solution of \eqref{eq:sys}, but rather it is a solution of a \mbox{finite-dimensional} approximation.

\begin{thm}[Tracking solution v-approximation] \label{th:main2} Let the assumptions of Theorem \ref{th:main} hold and consider the v-approximation of the discrete-time FOS \eqref{eq:sys} from Definition \ref{def:vApp}. Moreover, let
\begin{align}
\|u_r(k)\|\leq b_{u_r}, && \|x_r(k)\|\leq b_{x_r}, && k\in\mathbb{N}\label{eq:brefs}
\end{align}
for two scalars $b_{x_r},b_{u_r}\in\mathbb{R}^+$. Then, selecting $A=\tilde{A}_v$, $B=\tilde{B}_v$, and $C = [I_{n_x\times n_x},0_{n_x\times 1},\dots,0_{n_x\times 1}]$, the control law
\begin{align}
 u(k)=u_r(k)+K_v\tilde{e}(k),  && \tilde{e}(k)=\tilde{x}(k)-x_e(k)\label{eq:utrack}
\end{align}
solves Problem \ref{pb:2} with \begin{align}
d&= {\kappa\over 1-\kappa}\left(b_{x_r}\sum_{i=1}^l\|\hat{A}_{0}^{-1}A_i\| \phi_{{{a}}_i}(v) +b_{u_r}\sum_{i=1}^r \|\hat{A}_{0}^{-1}B_i\|\phi_{{{b}}_i}(v)\right) \label{eq:ddef}
\end{align}
and $\kappa$, $\gamma(\cdot)$, and $\phi_{{{a}}}(\cdot)$ as described in Theorem \ref{th:main}.
\end{thm}

It is worth noticing that the existence of bounded state and input trajectories $\{u_r(k)\}_{k\in\mathbb{N}}$ and $\{x_r(k)\}_{k\in\mathbb{N}}$ required in \eqref{eq:brefs} is always guaranteed by the Assumption~\ref{ass:controllability}. Moreover, we recall that since the exponential function can be rewritten as the infinite sum
$
e^x = \sum_{j=0}^\infty {x^j\over j!},
$
the function $\phi_{{a}}(v)$ defined in Theorem \ref{th:main} can be equivalently reformulated as 
\begin{align}
\phi_{{a}}(v)= e^{{a}} -\sum_{j=0}^{v} {{{a}}^j\over j!}=\sum_{j=v+1}^\infty {{{a}}^j\over j!} \label{eq:phialpha}
\end{align}
that decreases to zero as $v$ goes to infinity. As a consequence, the asymptotic bound of Theorem \ref{th:main2} is decreasing with the increment in the accuracy of the approximation. 

In practice, Theorem \ref{th:main2} shows that the state and input trajectories of the original FOS can track arbitrarily well (by increasing the value of $v$) any references $\{x_r(k)\}_{k\in\mathbb{N}}$ and $\{u_r(k)\}_{k\in\mathbb{N}}$ that are generated as a solution of \eqref{eq:exosys}, where $x_e$ is the state of the v-approximation. As a consequence, we can control the original infinite-dimensional FOS by designing a finite-dimensional controller (using standard techniques for finite dimensional systems) for the auxiliary system \eqref{eq:exosys}. The overall controller (that includes the v-approximation system) achieves in closed loop a tracking error whose ultimate bound decreases with the increment in the accuracy of the approximation. An illustrative example is presented in the following section with the design of a finite-dimensional controller that is able to steer the infinite-dimensional FOS around a desired trajectory.

\section{Numerical simulation} \label{s:numeric}
Consider the discrete-time \mbox{fractional-order} system \eqref{eq:sys} with $A = \begin{pmatrix}1&1\\0&1\end{pmatrix}$, $B_1 = \begin{pmatrix}0\\1\end{pmatrix}$, $l=3$, $A_1 = I_{2 \times 2}$, $A_2= A$, $A_3= -A$ , $a_1 = 0$, $a_2 = 1.7$, $a_3 = 0$, $r=1$, $b_1=0$, $s=1$, $g_1=0$, which satisfies Assumption \ref{ass:model} since $\sum_{i=1}^lA_i =I_{2\times 2} $ is a full rank matrix. The matrix $G_1$ is selected depending on the simulation.

\begin{figure*}
\begin{center}
\includegraphics[width=\textwidth]{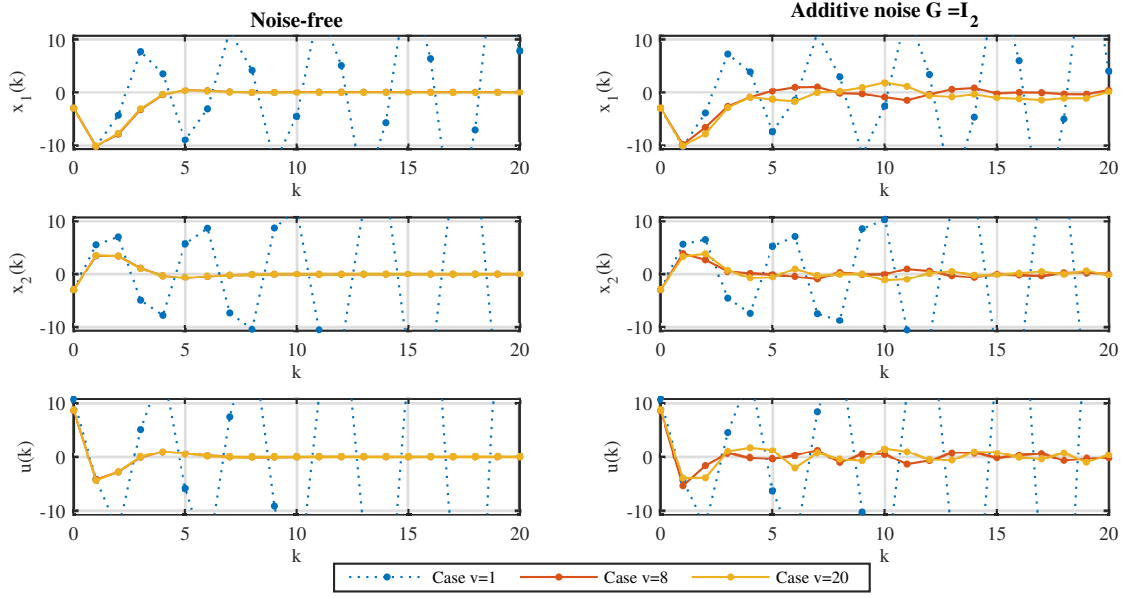}    
 
\caption{Closed-loop trajectories for the regulation of the origin in cases with and without noise (right and left column, respectively).\vspace{-.2cm}}  
\label{fig:1} 
\end{center}      
\end{figure*}

\begin{figure*}
\begin{center}
\includegraphics[width=\textwidth]{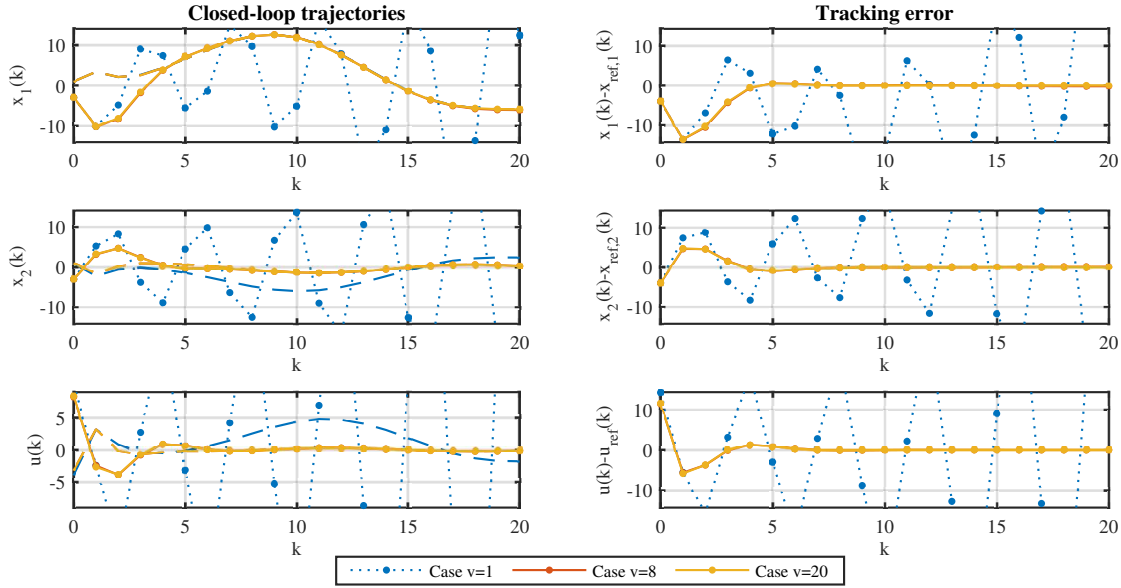}    
 
\caption{Closed-loop trajectories and reference signals (dashed lines) for the scenario of trajectory-tracking of the solution of an exogenous system. For the sake of clarity, the plot associated with the unstable case $v=1$ are dotted.\vspace{-.2cm}}  
\label{fig:2} 
\end{center}      
\end{figure*}

For different values of $v\in\mathbb{N}^{+}$, we display the behavior of the closed-loop \eqref{eq:sys} with \eqref{eq:input} and \eqref{eq:utrack} for both the cases of convergence to the origin and trajectory-tracking displayed in Fig. \ref{fig:1} and Fig. \ref{fig:2}, respectively. Moreover, in the first case, we consider both the noise-free case with $G_1=\mathbf{0}_{2\times 2}$ and the case with noise where $G_1=I_{2}$ and the noise $\{w(k)\}_{k\in\mathbb{N}}$ is generated as a random signal uniformly distributed within the interval $[-0.5,0.5]$ (i.e., $b_w=0.5$). 

Different approaches can be used to generate solutions of \eqref{eq:exosys}, that are reference trajectories for the trajectory-tracking algorithm. In this example, we used a model predictive control scheme (detailed in Appendix~\ref{app:mpc}) to compute reference trajectories for the $v$-approximation that drives $x_1$ toward a pre-defined sinusoidal trajectory. It is worth noticing that for each $v$, we generally have a different state and input solution of the $v$-approximation \eqref{eq:exosys} that achieve our objective (i.e., that drives $x_1$). This behavior is shown in Fig. \ref{fig:2} where the reference trajectory for $x_2$ in the case $v=1$ (dashed-line blue) differs from the case $v=9$ (dashed-line orange).  Despite this, in both cases the state $x_1$ is driven toward the desired pre-defined sinusoidal trajectory.

It is worth noticing that the system under analysis is open-loop unstable. Moreover, despite the $v$-approximation with $v=1$ being stabilizable, it follows that the associated linear controller does not stabilize the original infinite-dimensional system. Specifically, for the case $v=1$, Assumption \ref{ass:controllability} is satisfied but the condition \eqref{eq:conditionEps} is violated. Fig. \ref{fig:1} and Fig. \ref{fig:2} show that the closed-loop system is unstable with the error trajectory oscillating with increasing magnitude. To avoid this undesired behavior, Theorem \ref{th:main} provides a suitable selection of the value of $v$ that satisfies \eqref{eq:conditionEps} and therefore guarantees global convergence to the origin and, for the case of trajectory-tracking, global convergence to an ultimate bound $d$ around the reference signal.

\begin{figure}
\begin{center}
\includegraphics[width=0.5\textwidth]{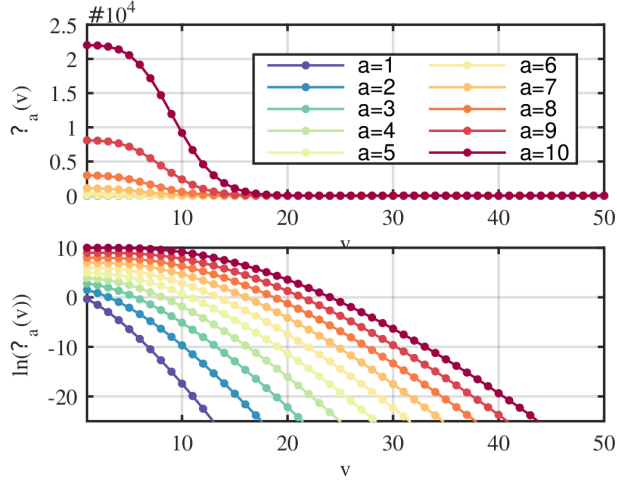}    
\caption{On the top we show the values of $\phi_a(v)$ by varying $v$ and $a$, and on the bottom we show these in a logarithmic scale.}
\label{fig:phia} 
\end{center}      
\end{figure}

It is worth noticing that the function $\phi_{{a}}(v)$ strongly decreases as the value of $v$ increases, as depicted in Fig. \ref{fig:phia}. This plays an important role in the design and analysis of the controller, since small values of $\phi_{{a}}(v)$ render condition \eqref{eq:conditionEps} easy to satisfy and lead to small values of the term $d$ in \eqref{eq:ddef}. This is reflected in the following table where the value of $d$ is close to zero since $v=8$. 

\begin{table}[h!]
	\centering
	\begin{tabular}{|l|l|l|l|}
		\hline
		$v$  & $1$ & $8$ & $20$ \\
		\hline
		$\Psi(v)$ & $4.4883$ & $6.3481\cdot 10^{-4}$ & $1.5\cdot 10^{-15}$ \\
		$c_\Psi$ & $163.2081$ & $647.1015$ & $1495.2155$ \\
		\hline
		$c_\gamma$  & $-$ & $462.2394$ & $15.1032$ \\
		$d$ & $-$ & $5.768\cdot 10^{-3}$ & $1.9\cdot 10^{-16}$ \\
		\hline
	\end{tabular}
\end{table}

\par Besides, the table highlights that the constants in Theorem \ref{th:main} are generally conservative. For instance, this applies to the term $c_\gamma$ that leads to a pessimistic estimate of the effect of the disturbances, which is in practice much smaller, as illustrated in Fig. \ref{fig:1}.

\section{Conclusions}
This work introduces a set of procedures to design \mbox{finite-dimensional} controllers for FOSs that guarantee closed-loop convergence to the origin, in the case of regulation, or to a reference trajectory, in the case of trajectory-tracking. This is achieved exploiting a linear controller designed for a suitable \mbox{finite-dimensional} approximation of the original FOS. For the case of trajectory-tracking, we consider both cases where the trajectories to be tracked are solutions of a FOS or solution of an exogenous \mbox{finite-dimensional} linear system. The latter case is desired in many practical cases where the reference trajectories are to be generated. All the proposed results are global, in the sense that the convergence is guaranteed for any initial condition of the system, and robust to additive disturbances.

\appendix
\section{Appendix} 
The following lemma is used in the proofs of the main results.

\begin{lem} \label{lem:bounLypa} Consider the system \eqref{eq:augsys} and let Assumptions \ref{ass:model} and \ref{ass:controllability} hold. Then, for any positive definite matrix $Q\succ 0$ there exists a pair of matrices $K_v$ and $P\succ 0$, with $P$ being a symmetric matrix, such that \eqref{eq:lyap} holds and, for the closed-loop system \eqref{eq:augsys} with \eqref{eq:input}, the state trajectory satisfies

\begin{align}
\|\tilde{x}(k_0+k)\|\leq \beta(\|\tilde{x}(k_0)\|,k) +c_\Psi \|r\|_{[k_0,k_0+k-1]}\label{eq:issmax}
\end{align}

with $ k_0,k\geq 0$, for a \mbox{class-$\mathcal{KL}$} function $\beta:\mathbb{R}^{+}\times \mathbb{N}^+\to\mathbb{R}^{+}$, and $c_\Psi$ from Theorem \ref{th:main}.
\end{lem}

\begin{pf}
Consider the function
\begin{align}
V(\tilde{x}) = \tilde{x}^\top P\tilde{x} \label{eq:Vdef}.
\end{align}
Using the system model \eqref{eq:augsys} and the equality \eqref{eq:lyap}, the function \eqref{eq:Vdef} evolves in time according to
\begin{align}
V^+ &= \tilde{x}^\top A_K^\top P A_K\tilde{x} +  r^\top \tilde{G}_v^\top P\tilde{G}_v r + 2 r^\top \tilde{G}_v^\top PA_K\tilde{x}\nonumber\\
&= V-(1-\theta)\tilde{x}^\top Q\tilde{x}-\theta \tilde{x}^\top Q\tilde{x} +  r^\top \tilde{G}_v^\top P\tilde{G}_v r + 2 r^\top \tilde{G}_v^\top PA_K\tilde{x}\label{eq:Vlb1}
\end{align}
for any scalar $\theta\in(0,1)$. Notice that for any $\bar{c}_\alpha,\bar{c}_\beta,\bar{c}_\gamma\in\mathbb{R}$ the inequality
\begin{align}
-\bar{c}_\alpha x^2\!\!+\!\bar{c}_\beta y^2\!\!+\!\bar{c}_\gamma xy\!=\! -(\bar{c}_ax\!-\!\bar{c}_by)^2\!\!+\!\bar{c}_c y^2\!\leq \bar{c}_c y^2 \label{eq:compSquars}
\end{align}
holds with $\bar{c}_a=\sqrt{\bar{c}_\alpha}$, $\bar{c}_b={\bar{c}_\gamma\over 2\bar{c}_a}$, $\bar{c}_c=\bar{c}_\beta+\bar{c}_b^2$. Then, we can bound the last part of the right hand side of \eqref{eq:Vlb1} as follows
$$ -\theta \tilde{x}^\top Q\tilde{x} +  r^\top \tilde{G}_v^\top P \tilde{G}_vr + 2 r^\top \tilde{G}_v^\top PA_K\tilde{x} \leq - \theta \lambda_{\mathrm{min}}(Q)\|\tilde{x}\|^2 + \|r\|^2\lambda_{\mathrm{max}}(\tilde{G}_v^\top P \tilde{G}_v) + \\2 \|r\| \|\tilde{G}_v^\top PA_K\| \|\tilde{x}\| \leq c_2\|r\|^2$$ with $c_2 = \lambda_{\mathrm{max}}(\tilde{G}_v^\top P \tilde{G}_v) +{\|\tilde{G}_v^\top PA_K\|^2\over \theta \lambda_{\mathrm{min}}(Q)}$ and where in the first inequality we used the fact that for any symmetric matrix $Q$ the inequality $$\lambda_{\min}(Q)\|x\|^2\leq x^{\top}Qx \leq \lambda_{\max}(Q)\|x\|^2$$ holds and in the last inequality, we used \eqref{eq:compSquars} with $\bar{c}_\alpha=\theta \lambda_{\mathrm{min}}(Q)$, $\bar{c}_\beta=\lambda_{\mathrm{max}}(\tilde{G}_v^\top P\tilde{G}_v)$ and $\bar{c}_\gamma = 2 \|\tilde{G}_v^\top PA_K\|$. As a consequence, by replacing it in \eqref{eq:Vlb1} we obtain $$ V(\tilde{x}(k+1))\leq V(\tilde{x}(k)) -(1-\theta)\tilde{x}(k)^\top Q\tilde{x}(k) + c_2\|r(k)\|^2 .$$
 Therefore, we can write $$\alpha_1(\|\tilde{x}\|)\leq V(\tilde{x})\leq \alpha_2(\|\tilde{x}\|),$$ $$V(\tilde{x}(k+1)) -V(\tilde{x}(k)) \leq  -\alpha_3(\|\tilde{x} \|) + \sigma(\|r(k)\|)  \leq  -\alpha_4(V(\tilde{x})) + \sigma(\|r(k)\|)$$ with $\alpha_1(r) = \lambda_{\mathrm{min}}(P)r^2$, $\alpha_2(r) = \lambda_{\mathrm{max}}(P)r^2$,  $\alpha_3(r)~=~(1-\theta)\lambda_{\mathrm{min}}(Q)r^2$, $\sigma(r)=c_2 r^2$, and $\alpha_4 = \alpha_3 \circ \alpha_2^{-1}$, i.e., $\alpha_4(r)=c_4r$, with $c_4={(1-\theta)\lambda_{\mathrm{min}}(Q) \over \lambda_{\mathrm{max}}(P)}$. Noticing that $\alpha_1,\alpha_2,\alpha_3,\sigma\in\mathcal{K}_\infty$, the function $V(\cdot)$ is an \mbox{ISS-Lyapunov} function (see, e.g., Definition 3.2 of \cite{Jiang2001}), which implies that there exist a \mbox{class-$\mathcal{KL}$} function $\beta:\mathbb{R}^+\times \mathbb{R}^+\to\mathbb{R}^+$ and a class-$\mathcal{K}$ function $\gamma:\mathbb{R}^+\to\mathbb{R}^+$ such that 
\begin{align}
\|\tilde{x}(k_0+k)\|\leq \beta(\|\tilde{x}(k_0)\|,k) + \gamma(\|r\|_{[k_0:k_0+k-1]}) \label{eq:issss}
\end{align}
for all $k_0, k\geq 0$. Now, let $\hat{\alpha}_4(r)$ be any class-$\mathcal{K}_\infty$ lower bound of $\alpha_4\in\mathcal{K}_\infty$ such that $\Id-\hat{\alpha}\in\mathcal{K}$, e.g., $\hat{\alpha}_4(r)~=~\hat{c}_4r$ with $\hat{c}_4=\min(c_4,\hat{\theta})$ for any $\hat{\theta}\in(0,1)$. Moreover, let $\rho(r)=c_\rho r$, with $c_\rho\in (0,1)$. Then, following \cite{Jiang2001} (specifically, from (13) to (17)) it is possible to show that \eqref{eq:issss} holds with $\gamma= \alpha_1^{-1}\circ\hat{\gamma}$  where $\hat{\gamma}=\hat{\alpha}_4^{-1}\circ \rho^{-1}\circ \sigma$, i.e., $\hat{\gamma}(r)=   {c_2\over \hat{c}_4 c_\rho} r^2$, and $\beta(s,t) = \alpha_1^{-1}(\hat{\beta}(\alpha_2(s),t)) $ for a \mbox{class-$\mathcal{KL}$} function $\hat{\beta}:\mathbb{R}^{+}\times \mathbb{R}^+\to\mathbb{R}^{+}$, which concludes the proof. \hfill $\blacksquare$
\end{pf}

\subsection{Proof of Theorem \ref{th:main}}
The proof is structured as follows. First, we show that the state trajectory of the closed-loop system \eqref{eq:augsys} with \eqref{eq:input} is uniformly bounded over time and converges to an ultimate bound proportional to the value of $\|r\|_\infty$, which also includes the contributes of the disturbance vectors $\{w(k)\}_{k\in\mathbb{N}}$. Then, noticing that the magnitude of $\|r(k)\|$ is increasing with the increase of the norm of the past state vectors, a contraction argument can be invoked to show convergence to an ultimate bound that is proportional only to $\|w\|_\infty$, which concludes the proof.

\subsubsection{Ultimate boundedness}

Taking the norm of  $r(k)$ in \eqref{eq:defd} and using the bound~\eqref{eq:bound_w} on the disturbance, results in $$ \|r(k)\|\leq  \sum_{i=1}^l\|\hat{A}_{0}^{-1}A_i\| \sum_{j=v+1}^{\infty} |c_j^{{{a}}_i}|\|x\|_{[0,k-v]}\\ + \sum_{i=1}^r \|\hat{A}_{0}^{-1}B_iK_v\|\sum_{j=v+1}^{\infty} |c_j^{{{b}}_i}|\|x\|_{[0,k-v-1]} + \gamma_w $$ with
\begin{align}
\gamma_w \!=\!b_w\!\sum_{i=1}^s \!\|\hat{A}_{0}^{-1}G_i\|\!\sum_{j=0}^\infty {{{g}}_i^j\over j!}\!=\! b_w\!\sum_{i=1}^s \|\hat{A}_{0}^{-1}G_i\| e^{{{g}}_i} .\label{eq:defgamma}
\end{align}
Notice that we used the fact that $\|x\|_{[-\infty,k-v]} = \|x\|_{[0,k-v]}$ since, as described in the problem formulation, the state of the system is considered to be the zero vector before the initial time $0$. Combining the latter with \eqref{eq:boundc} and the defining \eqref{eq:phialpha} we obtain
\begin{align}
&\|r(k)\|\!\leq\!  \Psi(v) |\|x\|_{[0,k-v]} \!+\! \gamma_w \leq \Psi(v) |\|x\|_\infty \!\!+\! \gamma_w\label{eq:iss1}
\end{align}
with $\Psi(\cdot)$ as in \eqref{eq:defPsi}.
Combining \eqref{eq:issmax} in Lemma \ref{lem:bounLypa} with \eqref{eq:iss1} results in
$$\|\tilde{x}\|_\infty\leq \beta(\|\tilde{x}_0\|,0)+ c_\Psi \|r\|_\infty \leq \beta(\|\tilde{x}_0\|,0)+ c_\Psi \Psi(v) |\|x\|_\infty + c_\Psi\gamma_w,$$ that under condition \eqref{eq:conditionEps} leads to
\begin{align}
\|\tilde{x}\|_\infty\!\leq  (1\!-\!c_\Psi \Psi(v) )^{-1}(\beta(\|\tilde{x}_0\|,0) \!+\!  c_\Psi \gamma_w \!)\!=:\!b_x .\label{eq:defbx}
\end{align}
Replacing the latter bound on the state trajectory into  \eqref{eq:iss1}, and then \eqref{eq:issmax}, results in 
\begin{align}
\|\tilde{x}(k)\|\!\leq\! \beta(\|\tilde{x}_0\|,k) \!+\! s_x, \;s_x:=c_\Psi [\Psi(v)b_x + \gamma_w ] \label{eq:bund0}
\end{align}
for a \mbox{class-$\mathcal{KL}$} function $\beta:\mathbb{R}^{+}\times \mathbb{R}^+\to\mathbb{R}^{+}$.\label{eq:convergence}

\subsubsection{Convergence}

We start by noticing that for any $M\in\mathbb{N}^+$ the bound on the norm of $r(k)$ in \eqref{eq:defd} can be decomposed as
\begin{align}
&\|r(k)\|\leq  \sum_{i=1}^l\|\hat{A}_{0}^{-1}A_i\| \sum_{j=v+1}^{v+1+M} |c_j^{{{a}}_i}|\|x\|_{[k-v-M,k-v]} \nonumber \\
&\quad+ \sum_{i=1}^l\|\hat{A}_{0}^{-1}A_i\| \sum_{j=v+2+M}^{\infty} |c_j^{{{a}}_i}|\|x\|_{[0,k-v-M-1]} \nonumber \\
&\quad+ \sum_{i=1}^r \|\hat{A}_{0}^{-1}B_iK_v\|\sum_{j=v+1}^{v+1+M} |c_j^{{{b}}_i}|\|x\|_{[k-v-M-1,k-v-1]}\nonumber\\
&\quad+\sum_{i=1}^r \|\hat{A}_{0}^{-1}B_iK_v\|\sum_{j=v+M+2}^{\infty} |c_j^{{{b}}_i}|\|x\|_{[0,k-v-M-2]} \nonumber\\
&\quad+ \gamma_w \nonumber\\
&\leq \!\Phi_v(M)\|x\|_{[k\!-\!v\!-\!M\!-\!1,k\!-\!v]} \!+\!\Psi(v\!+\!1\!+\!M)b_x \!+\!\gamma_w \label{eq:rboundec}
\end{align}
with $\gamma_w$ as in \eqref{eq:defgamma}, $b_x$ as in \eqref{eq:defbx}, $\Psi(\cdot)$ as in \eqref{eq:defPsi}, and
\begin{align}
\Phi_v(M):=&\sum_{i=1}^l\|\hat{A}_{0}^{-1}A_i\| \sum_{j=v+1}^{v+1+M} |c_j^{{{a}}_i}|   + \sum_{i=1}^r \|\hat{A}_{0}^{-1}B_iK_v\|\sum_{j=v+1}^{v+1+M} |c_j^{{{b}}_i}|.
\end{align}
The remainder of the proof is carried out by recursion. Specifically, for a generic $i\in\mathbb{N}$, consider a constant $s_i\geq 0$ and a function $\beta_i:\mathbb{R}^{+}\times\mathbb{R}^{+}\to\mathbb{R}^{+}$ such that
\begin{align}
\|\tilde{x}(k)\|\leq \beta_i(\|\tilde{x}_0\|,k) + s_i, \label{eq:siss_si}
\end{align}
holds, where for any given $r\in\mathbb{R}^{+}$, $\beta_i(r,s)$ is bounded over $s\in\mathbb{R}^{+}$  and satisfies $\beta_i(r,s)\to 0$ as $s\to\infty$. Moreover, consider a $\theta_i\in\mathbb{R}$ that satisfies 
\begin{align}
\theta_i\! \in\!\left(\!0, \min\left\{1,\!{\kappa-\epsilon\over \epsilon}s_i\right\}\right)\!\!\implies\!\! \epsilon(s_i+\theta_i)\leq\kappa s_i\label{eq:thetaii}
\end{align}
with $\kappa \in(\epsilon,1)$ and $\epsilon=c_\Psi \Psi(v)$, which satisfies $\epsilon\in(0,1)$ as a consequence of \eqref{eq:conditionEps}, and select a $\theta_{i+1}\in(0, \kappa\theta_i)$, which from \eqref{eq:thetaii} satisfies 
\begin{align}
0<\epsilon(s_i+\theta_{i+1})<\epsilon(s_i+\theta_{i})<\kappa s_i.\label{eq:thetaiii}
\end{align}
From \eqref{eq:siss_si}, there exists a $k_i=k_i(\theta_{i+1})$ such that 
$
k\geq k_i \implies \|x\|_{[k-v+1,k]}\leq \|\tilde{x}(k)\|\leq s_i+{\theta_{i+1}\over 2}$, which implies that for any $M_i\in\mathbb{N}$ we have
$$ k\geq k_i + M_i+2 =:\bar{k}_i \implies \|x\|_{[k-v-M_i-1,k-v]} \leq s_i+{\theta_{i+1}\over 2}.$$
 Hence, combining the latter with \eqref{eq:rboundec}, we obtain $$ \|r(k)\| \leq \Phi_v(M_i)\left(s_i+{\theta_{i+1}\over 2}\right) +\Psi(v+M_i+1)b_x +\gamma_w$$
for all $k\geq \bar{k}_i$. Therefore, by choosing $M_i$ to satisfy
\begin{align}
c_\Psi\Psi(v+M_i+1)b_x\leq \epsilon \theta_{i+1}/2 \label{eq:Mii},
\end{align}
from \eqref{eq:issmax} in Lemma \ref{lem:bounLypa}, we have
\begin{align}
&\|\tilde{x}(\bar{k}_i + k)\| \leq \beta(\|\tilde{x}(\bar{k}_i)\|,k) +s_{i+1}\label{eq:betaip1}
\end{align}
with 
\begin{align}
s_{i+1} &= c_\Psi \left( \Phi_v(M_i)\left(s_i\!+\!{\theta_{i+1}\over 2}\right) \!+\!\Psi(v\!+\!M_i\!+\!1)b_x\!+\!\gamma_w\right)\nonumber\\
& \leq c_\Psi \Phi_v(M)\left(\!s_i+{\theta_{i+1}\over 2}\!\right) + \epsilon {\theta_{i+1} \over 2}
+c_\Psi\gamma_w\nonumber\\
& \leq c_\Psi \Psi(v) \left(s_i+{\theta_{i+1}\over 2}\right) + \epsilon {\theta_{i+1} \over 2}+c_\Psi\gamma_w\nonumber\\
& = \epsilon \left(s_i+{\theta_{i+1}\over 2}\right) + \epsilon {\theta_{i+1} \over 2}+c_\Psi\gamma_w\nonumber\\
& \leq \epsilon (s_i+\theta_{i+1})+c_\Psi\gamma_w\leq \kappa s_i+c_\Psi\gamma_w,\nonumber
\end{align}
where the first inequality comes from  \eqref{eq:Mii}, the second from the fact that $\Psi(v)\geq \Phi_v(M)$ for any $v,M\in\mathbb{N}$, and the third from \eqref{eq:thetaiii}. Combining \eqref{eq:siss_si} with \eqref{eq:betaip1}, we obtain
\begin{align}
\|\tilde{x}(k)\|\!\leq \!\begin{cases}
\beta_i(\|\tilde{x}_0\|,k) + s_i,\hfill k\in[0,\bar{k}_i-1]\\
\min\{\beta_i(\|\tilde{x}_0\|,k) \!+ \!s_i, \nonumber\\
\quad \beta(\beta_i(\|\tilde{x}_0\|,\bar{k}_i) \!+\! s_i,k\!-\!\bar{k}_i) \!+\!\kappa s_i\!+\!c_\Psi\gamma_w\},\nonumber\\\hfill k\in[\bar{k}_i,+\infty]
 \end{cases}
\end{align}
and consequently
\begin{align}
\|\tilde{x}(k)\|\leq \beta_{i+1}(\|\tilde{x}_0\|,k) +\kappa s_i+c_\Psi\gamma_w\
\end{align}
with 
\begin{align}
\beta_{i+1}(r,k)\!:=\! \begin{cases}
\beta_i(r,k) \!+\! (1\!-\!\kappa)s_i\!-\!c_\Psi\gamma_w,\hfill k\in[0,\bar{k}_i-1]\\
\min\{\beta_i(r,k) + (1-\kappa)s_i-c_\Psi\gamma_w, \nonumber\\
\quad \beta(\|\tilde{x}(\bar{k}_i)\|,k-\bar{k}_i) \},\hfill k\in[\bar{k}_i,+\infty]
 \end{cases}
\end{align}
that recursively satisfy the fact that $\beta_{i+1}(r,k)\to 0$ as $k\to\infty$. Roughly speaking, we showed that as the state converges toward $s_i$, it is possible to compute another ultimate bound $s_{i+1}\leq \kappa s_i+c_\Psi\gamma_w$. Iterating the process starting from \eqref{eq:bund0}, i.e., $\beta_0(r,s)=\beta(r,s)$, $s_0=s_x$, and any $ \theta_0 \in\left(0, \min\left\{1,{\kappa-\epsilon\over \epsilon}s_x\right\}\right)$,
we obtain 
$$s_i = \kappa^{i}s_0 + c_\Psi\gamma_w\left(\sum_{j=0}^{i-1}\kappa^j\right)$$
and therefore
\begin{align}
\|\tilde{x}(k)\|\leq \beta_\infty(\|\tilde{x}_0\|,k) +s_\infty \label{eq:sinfbound}
\end{align}
with 
$$
s_\infty = \lim_{i\to\infty} s_i = c_\Psi\gamma_w{\kappa \over 1-\kappa} \nonumber
$$
where $\beta_\infty(r,k)\to 0$ as $k\to\infty$. The proof is concluded by noticing that $\|\tilde{x}_0\|=\|x_0\|$ and $\|x(k)\|\leq\|\tilde{x}(k)\|$, which implies \eqref{eq:issdef} with $\beta = \beta_\infty$. \hfill $\blacksquare$

\subsection{Proof of Corollary \ref{cor:1}}

By the linearity of the Gr\"{u}nwald-Letnikov \mbox{fractional-order} difference operator, the evolution of the error vector $e(k)$ is governed by the \mbox{fractional-order} system 
$$ \sum_{i=1}^{l}A_i\Delta^{{{a}}_i}e(k+1) = \sum_{i=1}^{r}B_i\Delta^{{{b}}_i}(u(k)-u_r(k)) +  \sum_{i=1}^{s}G_i\Delta^{{{g}}_i}w(k).$$
Hence, the result readily follows by applying  Theorem \ref{th:main} with suitable change of state and input coordinates.
\hfill $\blacksquare$

\subsection{Proof of Theorem \ref{th:main2}}

For the closed-loop system \eqref{eq:sys} with \eqref{eq:utrack} the tracking error trajectory satisfies $$ \tilde{e}(k+1) = \tilde{A}_v \tilde{e}(x)+\tilde{B}_v(u(k)-u_r(k)) + \tilde{G}_v r(k)$$ where the norm of $r(k)$ is bounded as 
\begin{align}
&\|r(k)\|\leq \sum_{i=1}^l\|\hat{A}_{0}^{-1}A_i\| \sum_{j=v+1}^{\infty} |c_j^{{{a}}_i}|\|e(k-j+1)\|\nonumber \\
&\quad + \sum_{i=1}^l\|\hat{A}_{0}^{-1}A_i\| \sum_{j=v+1}^{\infty} |c_j^{{{a}}_i}|\|x_r(k-j+1)\|\nonumber \\
&\quad+\sum_{i=1}^r \|\hat{A}_{0}^{-1}B_iK_v\|\sum_{j=v+1}^{\infty}|c_j^{{{b}}_i}|\|e(k-j)\|
\nonumber \\
&\quad+\sum_{i=1}^r \|\hat{A}_{0}^{-1}B_i\|\sum_{j=v+1}^{\infty}|c_j^{{{b}}_i}|\|u_r(k-j)\|+\gamma_w\label{eq:defw2}
\end{align}
In contrast to the case considered in Corollary \ref{cor:1}, where the reference to track is a solution of the FOS \eqref{eq:sys}, the term $r(k)$ it is not only function of the disturbance $\{w(k)\}_{k\in\mathbb{N}}$ and the previous error vectors $\{e(k)\}_{k\in\mathbb{N}}$ but also of the reference signals $\{x_r(k)\}_{k\in\mathbb{N}}$ and $\{u_r(k)\}_{k\in\mathbb{N}}$. The same proof of Theorem \ref{th:main} applies here, the only difference being the term $\gamma_w$, which is to be replaced with 
$$\sum_{i=1}^l\|\hat{A}_{0}^{-1}A_i\| \sum_{j=v+1}^{\infty} |c_j^{{{a}}_i}|\|x_r(k-j+1)\|+\sum_{i=1}^r \|\hat{A}_{0}^{-1}B_i\|\sum_{j=v+1}^{\infty}|c_j^{{{b}}_i}|\|u_r(k-j)\|+\gamma_w$$ .
Using the latter substitution in \eqref{eq:sinfbound} evaluated on the error space results in 
\begin{align}
\|\tilde{e}(k)\|\leq \beta_\infty(\|\tilde{e}_0\|,k) + {\kappa\over 1-\kappa}\gamma_w + d
\end{align}
with $d$ in \eqref{eq:ddef} where we used the bounds \eqref{eq:brefs} and \eqref{eq:boundc} and the function $\phi_{{{a}}}(\cdot)$ from \eqref{eq:phialpha}. The proof is closed noticing  that $\|\tilde{e}(0)\|=\|e(0)\|$ and $\|e(k)\|\leq \|\tilde{e}(k)\|$ and therefore \eqref{eq:iips} holds with $\beta=\beta_\infty$ and $\gamma(\cdot)$ from Theorem \ref{th:main}.
\hfill $\blacksquare$

\section{Generation of reference trajectories} \label{app:mpc}

In practical applications, the trajectories $\{x_r(k)\}_{k\in\mathbb{N}}$ and $\{u_r(k)\}_{k\in\mathbb{N}}$ that are a solution of \eqref{eq:exosys} are unknown, and need to be computed. A possible approach consists of using model predictive control (MPC) schemes. For instance, consider the case where we wish to drive the first component of the state of the FOS system \eqref{eq:sys} to a given desired trajectory  
$$p_d(k) = -10\sin(0.2k)+3\sin(0.5 k), $$
possibly not feasible for \eqref{eq:exosys}, while minimizing the use of the input. Then, we can design a MPC scheme that minimizes the stage cost $$l_e(k,x,u) = 10\|x_1(k)-\!p_d(k)\|^2 \!+\! \|u(k)\|^2$$ where $x(k)=[x_1(k),x_2(k),\dots]'$. Specifically, for a given a pair $(k,x_e)\!\in\!\mathbb{R}_{\geq k_0} \times \mathbb{R}^{n_{x_e}}$ and an integer horizon length $N>0$, the open-loop optimization problem $\mathcal{P}(k,x_e)$ consists of finding the optimal control trajectory $$\bbaru_r^*\!:=\!\{u_r^*(k),u_r^*(k\!+\!1),\dots,u_r^*(k\!+\!N\!-\!1)\}$$ that solves
	\begin{align} 
	J^*(k,x_e) &= \min_{\boldsymbol{\bar{u}}_r}   J(k,x_e,\boldsymbol{\bar{u}}_r) \nonumber\\
	\text{s.t.} \quad & \bar{x}_e(i\!+\!1)\!=\!\tilde{A}_v\bar{x}_e(k)\!+\!\tilde{B}_v\bar{u}_r(i),&  i \in \mathbb{Z}_{k:k+N-1}\label{eq:sysAref}\\
	\quad& \bar{x}_e(k) = x_e,\nonumber
	\end{align}
	with $\bbaru_r=\{\bar{u}_r(k),\bar{u}_r(k+1),\dots,\bar{u}_r(k+N-1)\}$ and where the performance index is defined as
	\begin{align}
	J(k,\!x_e,\!\boldsymbol{\bar{u}}_r)\!:= \!\!\!\!\!\! \sum_{i=k}^{k+N-1} \!\!\!\!l(i,\!\bar{x}_e(i),\!\bar{u}_r(i))\!+\! m_s(\bar{x}_e(k\!+\!N)\!)\label{eq:costsJr}
	\end{align}
	with $l(k,x,u) =  c_{l_e}\atan (l_e(i,x,u)/c_{l_e}) + l_s(x)$, $m_s(x) = c_sx^{\top}Px$, and $l_s(x)=c_sx^{\top}Qx$ for some constants $c_{l_e},c_s>0$, and the matrices $P$ and $Q$ satisfy Assumption \ref{ass:controllability}. It is worth noticing that for small values of $c_s$ and high values of $c_{l_e}$, the effect of the function $l_s$ and $m_s$ and the saturation $c_{l_e}\atan (\cdot/c_{l_e})$ is neglectable for contained values of $\|x_e\|$. These terms are added to provide a formal guarantee of boundedness of the closed-loop trajectory, by following a similar reasoning to that presented in  \cite{AlessandrettiAguiarJonesTAC17} for continuous-time systems. Indeed, for bounded values of the state, arbitrarily small values $c_s$ and arbitrarily high values of $c_{l_e}$, the performance index is approximately given by
	\begin{align}
	J(k,x_e,\boldsymbol{\bar{u}}_r) &\sim  \sum_{i=k}^{k+N-1} l_e(i,\bar{x}_e(i),\bar{u}_r(i)).
	\end{align} 
Notice that, in the numerical example in Section \ref{s:numeric}, we used $c_{l_e}=1\cdot 10^5$ and $c_s=1\cdot 10^{-5}$.
	
	In order to make explicit the dependence of the optimal solution to the parameters of the open-loop optimization problem, let $\bar{u}^{*,k}(i)$ and $\bar{x}_e^{*,k}(i)$ denote the optimal input and state trajectories $\bbaru_r^{*,k}$ and $\bbarx_e^{*,k}$ computed by solving $\mathcal{P}(k,x_e(k))$ associated with the predicted time $i$. Then, the state and input reference trajectories $\bx_r$ and $\bu_r$ are computed as the output and input trajectories, respectively, of the closed-loop \eqref{eq:exosys} with
\begin{align} 
u_r(k)=\kappa_{MPC}(k,x_e(k)) := \bar{u}_r^{*,k}(k) \label{eq:MPC}
\end{align}
obtained by solving the open-loop optimization problem at any time step and applying the first optimal input to the system.

\subsection{Boundedness of the reference trajectories}

In contrast with the work \cite{AlessandrettiAguiarJonesTAC17}, where boundedness was obtained in continuous-time systems, in this section we provide the proof for the adopted discrete-time setting. Specifically, we show that the state and input reference signals associated with the closed-loop \eqref{eq:MPC} with \eqref{eq:sysAref} are bounded. For the sake of clarity, wherever clear from the context, we use the compact notation, e.g., $l(i,\bar{x}_e^{*, k}(i),\bar{u}_r^{*, k}(i)) = \bar{l}^{*, k}(i)$.

Consider the MPC value function 
\begin{align}
V_{MPC}(k,x_e) :=J(k,x_e,\bbaru_r^{*,k})\label{eq:valuef1}
\end{align}
and, for any $i\geq k$, we denote by the extended input trajectory the concatenation of the optimal input trajectory $\bbaru_r^{*,k}$ with the linear controller $u_r(k)=K_v x_e(i)$, that is 
\begin{align}
u^{ext,k}(i) &= \begin{cases}
\bar{u}_r^{*,k}(i),&\!\!i=k,\dots,k+N-1,\\
K_v x_e^{ext,k} (i), &\!\!i\geq k +N ,
\end{cases}\label{eq:ext}
\end{align}
where $x^{ext,k}(i)$ denotes the associated state trajectory.
\par \textbf{Value function evolution.} Multiplying the first equality in \eqref{eq:lyap} by $c_s$ and using the definition \eqref{eq:ext}, for any $i\geq k+N$ the following holds
\begin{align}
m_s^{ext,k}(i)= l_s^{ext,k}(i) + m_s^{ext,k}(i+1). \label{eq:mlineq}
\end{align}
Then, the MPC value function evolves as follows:
\begin{align}
&V_{MPC}(k+1,\;x_e(k+1))\leq J(k+1,\;x_e(k+1),\bu^{ext, k})\nonumber\\
&= \sum_{i=k+1}^{k+N} l^{ext, k}(i) + m^{ext, k}(k+N+1) +V_{MPC}(k,\;x_e(k))-V_{MPC}(k,\;x_e(k))\nonumber\\
&=V_{MPC}(k,\;x_e(k))-l^{ext, k}(k) -m_s^{ext, k}(k+N)+ m_s^{ext, k}(k+N+1) +l^{ext, k}(k+N)\nonumber\\
&= V_{MPC}(k,\;x_e(k))-l_s^{ext, k}(k) +c_{l_e}\atan (l^{ext, k}_e(k+N)/c_{l_e}) -c_{l_e}\atan (l^{ext, k}_e(k)/c_{l_e})\nonumber\\
&\leq V_{MPC}(k,\;x_e(k))-\alpha(\|x_e(k)\|) + 2c_{l_e}{\pi\over 2}\label{eq:Lyapdec}
\end{align}
where the first inequality is a consequence of the from the sub-optimality of $\bu^{ext, k}$, the last equality from \eqref{eq:mlineq}, the last inequality from the bounds
\begin{align}  
\|c_{l_e}\atan (l_e(i,x,u)/c_{l_e}) \|\leq c_{l_e}{\pi\over 2},\label{eq:boundatan}
\end{align}
and $x_e^{\top}Qx_e\geq \alpha(\|x_e\|)$ with $\alpha(r)=:\lambda_{\mathrm{min}}(Q)r^2$.
\par \textbf{Upper bound value function.}
Let $\bx^{K_v}$ and $\bu^{K_v}$ be the state and input trajectories associated with the system \eqref{eq:sysAref} in closed-loop with $u_r(k)=K_v x_e(k)$. Then, we can write
$$
m_s^{K_v}(i)= l_s^{K_v}(i) + m_s^{K_v}(i+1)
$$ that results in 
$$
\sum_{i=k}^{k+N-1} l_s^{K_v}(i) +  m_s^{K_v}(k+N) = m_s^{K_v}(k)\leq \alpha_c(\|x_e\|)$$, with $\alpha_c(r):=\lambda_{\mathrm{max}}(P) r^2$. By the sub-optimality of $\bu^{K_v}$, and combining the latter equality with the bound \eqref{eq:boundatan} and the definition \eqref{eq:costsJr}, it follows that to 
\begin{align} 
&V_{MPC}(k,\;x_e)=J(k,x_e,\bu^{*,k}) \leq J(k,x_e,\bu^{K_v})\leq \alpha_c(\|x_e\|)+Nc_{l_e}{\pi\over 2}.  \label{eq:ubVmpc}
\end{align}
\par \textbf{Lower bound value function.} From the definition \eqref{eq:costsJr} and the fact that $l_e(\cdot)$ is always non negative, we have
\begin{align}
V_{MPC}(k,x_e(k))=J(k,x_e(k),\bu^{*,k})\geq \alpha(\|x_e\|). \label{eq:lbVmpc}
\end{align}
\textbf{Shifted value function.}
Let us consider the shifted value function $$W(k,x_e):= V_{MPC}(k,x_e)- Nc_{l_e}{\pi\over 2}.$$ Then, using the bounds \eqref{eq:ubVmpc} and \eqref{eq:lbVmpc}, we can write  
$$\alpha(\|\hat{x}\|)-c_{l_e}{\pi\over 2}N\leq W(k,\;x_e(k)))\leq \alpha_c(\|x_e\|).$$
Moreover, from \eqref{eq:Lyapdec}, we obtain 
 $$W(k+1,\;x(k+1)))
 \leq W(k,\;x_e(k)) -\alpha(\|x_e(k)\|) + 2c_{l_e}{\pi\over 2}\leq W(k,\;x_e(k)) -\alpha(\alpha_c^{-1}(W(k,\;x_e(k)))) + 2c_{l_e}{\pi\over 2}$$. At this point, following standard arguments from input-to-state stability analysis for discrete time systems \cite{Jiang2001}, it is possible to show that $x_e$ is bounded, which implies that $x_r$ and $u_r$ are bounded. Specifically, let $\alpha_4$ be a class $\mathcal{K}_\infty$ function such that $\alpha_4(r)\leq\alpha(\alpha_c^{-1}(r))$, $\Id-\alpha_4$ belongs to class $\mathcal{K}_\infty$, and $\rho$ be a class-$\mathcal{K}_\infty$ function such that $\Id-\rho\in\mathcal{K}_\infty$ . Then, we obtain that 
\begin{align}
W(k)\leq \max(\beta_e(W(0),k),\hat{\gamma}(\alpha_3(N,n)))
\end{align}
with $\hat{\gamma}(r)=\alpha_4^{-1}(\rho^{-1}(r))$ and a class-$\mathcal{KL}$ function $\beta_e$. More on how to compute $\beta_e$ can be found in \cite{Jiang2001}.

\bibliographystyle{plain}        

\bibliography{FractionalOrder_v17} 

\end{document}